\documentclass[12pt]{article}
\usepackage{amsmath,amsfonts,amssymb}
\usepackage{secdot}
\usepackage{amsthm}
\usepackage{color}
\newtheorem{theorem}{\hspace*{\parindent}Theorem}

\newcounter{theremark}

\title{Generalized reduced modulus in Euclidean space}
\author{A.S. Afanaseva-Grigoreva$^{\rm a,b}$\footnote{Corresponding author. E-mail: A.S. Afanaseva-Grigoreva -- \emph{a.s.afanasevagrigoreva@yandex.ru}, K.A. Gulyaeva-- \emph{ karinagulyaeva8@gmail.com},
E.\:Prilepkina --  \emph{pril-elena@yandex.ru}}, K.A. Gulyaeva$^{\rm a,b}$, 
E.G.\:Prilepkina$^{\rm a}$
\\[5pt]
\\
\small{\textit{$\phantom{1}^a$Institute of Applied Mathematics,
FEBRAS, 7 Radio Street, Vladivostok,  690041, Russia}}
\\
\small{\textit{$\phantom{1}^b$Far Eastern Federal University, Ajax Bay~10, 690922 Vladivostok, Russia}}
}
\date{}
\begin{document}
\maketitle


\begin{abstract}
In Euclidean space condensers with variable potential levels and the presence of a free part at the boundary are studied. The asymptotic formula of the modulus of such condenser is obtained when the plates are pulled into points. The generalized reduced modulus is defined as the second term of this asymptotics. The relation between the reduced modulus and the discrete energy of the Robin function, monotonicity, and polarization behavior are established.   
\end{abstract}

\bigskip


\bigskip


\section{Introduction}

The notion of reduced modulus is well-known in geometric function theory. G.V.Kuzmina, E.G.Yemelyanov, and A.Yu.Solynin have been systematically developing extremal-metric approach to the definition of reduced modulus \cite{Kuz}. We hold on to the capacitative interpretation, which implies that the reduced modulus appears in the asymptotics of the condenser capacity while pulling its plates into points. The most complete description of the capacitative approach is presented in the research papers by V.N.Dubinin and his students (e.g., \cite{DubEir}). This paper can be viewed as a continuation of our collaborative research with V.N.Dubinin (\cite{DubPril8}, \cite{DubPril5}, \cite{DubPrilGreen}). In \cite{Dub3} the asymptotics of the condenser capacity with variable potential levels on a plane is studied. Now we study the analogue of such condenser in Euclidean space generalizing the asymptotics formula of \cite[p. 101]{DubPril8} by introducing variable potential levels and assuming the existence of a specialized set allocated at the domain boundary. It should be noted that the reduced modulus in Euclidean space was introduced by B.E.Levitsky \cite{Levic}. The interest for the studies of the reduced modulus stems from a variety of applications. In Euclidean space our results can be used in the studies of harmonic radius properties (\cite{Band}) and different discrete energies (\cite{Brauchart}, 
\cite{Borodachov}, \cite{Brauchart2}, \cite{Dubinin1},\cite{PrilAfan}). The fundamental properties of the reduced moduli are monotonicity, composition principles, and symmetrization transformation behaviors (\cite{DubEir}). In most cases the reduced modulus inherits these properties from the condenser capacity properties. Several symmetrization transformations for condensers with two plates are known so far. For condensers with three or more plates the symmetrization technique application faces severe challenges. In the current paper we construct the polarization of the generalized reduced modulus using the Dirichlet's integral of what we call the potential function for the reduced modulus expression. This approach traces back to \cite{DubPril5}, and we have extended it to Euclidean space case in \cite{GKP}.

\section{Generalized reduced modulus in Euclidean space}

Henceforth, $\mathbb{R}^n$ denotes $n$-dimensional Euclidean space comprised of points $\mathbf x=(x_1,\ldots,x_n)$, where $n\geq3$, $\omega_{n-1}$ is the area of a unitary hypersphere, $E(\mathbf a,r)=\{|\mathbf x-\mathbf a|\leq r\},$ $S(\mathbf a,r)=\partial E(\mathbf a,r)$. Let $B$ be a domain in $\mathbb{R}^n$. We call the generalized condenser in $\overline{B}$ the ordered collection 
$C=(\{F_0,F_1,\ldots,F_m \},\{\delta_0,\delta_1,\ldots,\delta_m \},B)$, where $F_k\subset\overline{B}$ -- closed nonempty nonoverlapping sets, $\delta_k$ -- real numbers, where $k = 0,1,\ldots,m,$ $ m\geq 1$. $F_k$ sets are called condenser $C$ plates, and $\delta_k$ are $F_k$ plate potential levels, where $k = 0,1,\ldots,m$. The capacity $\mathrm{cap}\,C$ of condenser $C$ is Dirichlet's integral infinum 
\begin{equation}
I(v,B):=\int_B|\nabla v|^2d\mathbf x
\end{equation}
among all functions $v$, which are continuous up to the $\partial {B}$, satisfy the Lipschitz condition locally in $B$, and equal to $\delta_k$ at $F_k$, where $k = 0,1,\ldots,m,d\mathbf x=dx_1\ldots dx_n$ (such functions $v$ are called \emph{admissible}). Condenser modulus $|C|$ is a quantity reciprocal to capacity: $|C|=(\mathrm{cap}\,C)^{-1}.$

It follows from the capacity definition that  we have the inequality
\begin{equation}\label{eq:monoton}
   \mathrm{cap}\,C\leq \mathrm{cap}\,\tilde C
\end{equation}
for condensers $C=(\{F_0,F_1,\ldots,F_m \},\{\delta_0,\delta_1,\ldots,\delta_m \},B)$ and \\ 
$\tilde C=(\{\tilde F_0,\tilde F_1,\ldots,\tilde F_m \},$ $\{\delta_0,\delta_1,\ldots,\delta_m \},B)$, where  $F_k \subset \tilde F_k$,  $k=0,\ldots,m.$ 

Let us define the condenser (considering $r$ is sufficiently small) 
$$C(r;B,X,\Delta,\Psi) =( \{\partial B,E(\mathbf x_1,\mu_1 r),\ldots,E(\mathbf x_m,\mu_m r)\}, \{0,\delta_1,\ldots,\delta_m\}),$$
 where $X=\{\mathbf x_k\}_{k=1}^m,\Delta=\{\delta_k\}_{k=1}^m,\ \Psi=\{\mu_k \}_{k=1}^m,E(\mathbf x_k,\mu_k r)=\{\mathbf x \in \mathbb{R}^n:|\mathbf x_k-\mathbf x|\leq\mu_k r\},$  $\delta_k\neq 0$, $k = 1,\ldots,m.$ In \cite{DubPril8} the reduced modulus of the domain $B$ with respect to the collections $X,\Delta$, and $\Psi$ is the limit 
 \begin{equation}\label{eq:redmod}
 M(B,X,\Delta,\Psi) = \lim\limits_{r \to 0}(|C(r;B,X,\Delta,\Psi)|-\nu\lambda_n r^{2-n}),
 \end{equation}
where $\nu=(\sum\delta _l^2\mu_l^{n-2})^{-1},$  $\lambda_n=1/((n-2)\omega_{n-1})$. Let us give more generalized definition of the reduced modulus in space examining the condenser with variable potential levels and substituting a "nullable plate" \ $\partial B$ with the plate $\Gamma\subset \partial B$. Let $\Gamma$ be a closed nonempty subset of $\partial B$, $X=\{\mathbf x_k\}_{k=1}^m$ --- different points of the domain $B$, $\Delta=\{\delta_k\}_{k=1}^m$--- numbers other than zero, $\Omega=\{\alpha_k\}_{k=1}^m$--- arbitrary numbers, $ \Psi=\{\mu_k\}_{k=1}^m$ --- positive numbers, $r>0$ is a sufficiently small number, 
\begin{equation}
\delta_k(r)=\delta_k+\alpha_kr^{n-2}+o(r^{n-2}),\ r \to 0,
\end{equation}
and closed sets $F_k(r)$ contain balls  and are contained in balls of the form
\begin{equation}\label{eq:cap}
 |\mathbf x-\mathbf x_k|\leq \mu_k r(1+O(r^{n-1})),r \to 0,
\end{equation}
$\mu_k>0$, where $k=1,\ldots,m$,  
\begin{equation}\label{eq:cap2}
C(r)=C(r;B,\Gamma, X,\Delta,\Omega,\Psi) =( \{\Gamma, F_1(r),\ldots,F_m(r)\}, \{0,\delta_1(r),\ldots,\delta_m(r)\}).
\end{equation}
In accordance with \eqref{eq:redmod}, the reduced modulus of the domain $B$ with regard to the collections $X,\Delta,$  $ \Omega $, $\Psi$ is the following limit
\begin{equation}\label{eq:mod2}
M(B,\Gamma, X,\Delta,\Omega,\Psi) = \lim\limits_{r \to 0}(|C(r;B,\Gamma, X,\Delta,\Omega,\Psi)|-\nu\lambda_n r^{2-n}).\end{equation}

  Let $B\subset \mathbb{R}^n$ be a bounded finitely connected domain with a piecewise smooth boundary, ${\mathbf x}_0\in B,$ nonempty set $\Gamma\subset \partial B$ consists of a finite number of closed piecewise smooth surfaces. Robin function  $g_B(\mathbf x,\mathbf x_0,\Gamma)$  with a pole $\mathbf x_0\in B$ is defined in \cite{GKP} as the function that is continuous in $\overline{B}\setminus \{\mathbf x_0\}$, continuously differentiable in  $\overline{B}\setminus (\Gamma\cup\{\mathbf x_0\})$, vanishes at $\Gamma,$ and at the remainder of the boundary $\partial B\setminus \Gamma$ its normal derivative equals zero. Moreover,  function 
 $$h_B(\mathbf x,\mathbf x_0,\Gamma)=g_B(\mathbf x,\mathbf x_0,\Gamma)-\lambda_n|\mathbf x-\mathbf x_0|^{2-n}$$
is harmonic in $B$.  
 The quantity $$r(B,\Gamma,\mathbf x_0)=((2-n)\omega_{n-1}h_B(\mathbf x_0,\mathbf x_0,\Gamma))^{1/(2-n)}$$ is called the Robin radius. If $\Gamma=\partial B$, the Robin function is identical to the Green function \cite{DubPril8}, and the Robin radius is identical to the harmonic radius \cite{Band}. Directly from the definition above, the following asymptotic expansion is derived 
\begin{equation}
    g_B(\mathbf x,\mathbf x_0,\Gamma)=\lambda_n(|\mathbf x-\mathbf x_0|^{2-n}-r(B,\Gamma,\mathbf x_0)^{2-n})+O(|\mathbf x-\mathbf x_0|),\ \mathbf x\to \mathbf x_0.
\end{equation}

\begin{theorem} Let $B\subset \mathbb{R}^n$ be a bounded finitely connected domain with a piecewise smooth boundary, $\Gamma\subset \partial B$ consists of a finite number of closed piecewise smooth surfaces. Then the reduced modulus \eqref{eq:mod2} can be expressed as 
$$M(B,\Gamma,X,\Delta,\Omega,\Psi)=\lambda_n\nu^2R(B,\Gamma,X,\Delta,\Omega,\Psi),$$ where 
\begin{multline}\label{basic}
R(B,\Gamma,X,\Delta,\Omega,\Psi)=-\sum\limits_{k=1}^m (\delta_k^2\mu_k^{2n-4}(r(B,\Gamma,\mathbf x_k))^{2-n}+2\alpha_k\delta_k\mu_k^{n-2})\\+ \sum\limits_{k=1}^m\sum\limits_{{l=1,} \atop{l\not=k}}^m \delta_k\delta_l\mu_k^{n-2}\mu_l^{n-2}(n-2)\omega_{n-1}g_B(\mathbf x_k,\mathbf x_l,\Gamma).\end{multline}
For the condenser capacity \eqref{eq:cap2} at  $r\to 0$ the following asymptotic formula takes place
\begin{equation}\label{ascon}
\frac{cap\,C(r)}{(n-2)\omega_{n-1}} = 
\left(\sum\limits_{k=1}^m\delta_k^2\mu_k^{n-2}\right)r^{n-2}-R(B,\Gamma,X,\Delta,\Omega,\Psi)r^{2n-4}+o(r^{2n-4}).
\end{equation}
\end{theorem}
\textbf{Proof.} 
Let us denote $\Psi_k(r) = \mu_k r$, $g(\mathbf x,\mathbf x_k)=(n-2)\omega_{n-1}g_B(\mathbf x,\mathbf x_k,\Gamma)$, $R_k=r(B,\Gamma,\mathbf x_k),$ where $\\ k=1,\ldots,m.$
Considering $r$ has a sufficiently small positive value, let us introduce an intermediate function
$$g(\mathbf x)=\sum\limits_{k=1}^m \delta_k(r) \sum\limits_{l=1}^m \beta_{lk}g(\mathbf x,\mathbf x_l),$$
where
$$
\beta_{lk}=
\begin{cases}
-\Psi_k^{n-2}(r)\Psi_l^{n-2}(r)g(\mathbf x_k,\mathbf x_l),& k\not=l,\\
\Psi_k^{n-2}(r)(1+\Psi_k^{n-2}(r)R_k^{2-n}),& k=l,
\end{cases}
$$
$l,k=1,\ldots,m.$ Function $g(\mathbf x)$ equals zero at $\Gamma$, and it has a singularity in each point $\mathbf x_k$, $k=1,\ldots,m.$ Now we want to consider the behavior of this function in the neighborhood of the point $\mathbf x_k$ for each $k=1,\ldots,m.$ For the sake of  clarity, let us expand $g(\mathbf x)$ in the following form
\begin{multline}
g(\mathbf x)=\delta_k(r)\beta_{kk}g(\mathbf x,\mathbf x_k)+\delta_k(r)\sum \limits_{{l=1,} \atop{l\not=k}}^m\beta_{lk}g(\mathbf x,\mathbf x_l)+\sum \limits_{{j=1,} \atop {j\not=k}}^m\delta_j(r)\beta_{kj}g(\mathbf x,\mathbf x_k)\\
+\sum \limits_{{j=1,} \atop {j\not=k}}^m\delta_j(r)\beta_{jj}g(\mathbf x,\mathbf x_j)+\sum \limits_{{j=1} \atop {j\not=k}}^m\delta_j(r)\sum \limits_{{l=1,l\not=j,} \atop {l\not=k}}^m\beta_{lj}g(\mathbf x,\mathbf x_l)
\end{multline}

Let us examine the sets $E_k=\{\mathbf x: 0,5\Psi_k(r)\leq |\mathbf x-\mathbf x_k|\leq 2\Psi_k(r), \ g(\mathbf x)/\delta_k (r) \geq 1\}.$  In the points  $\mathbf x\in  E_k,$ where $g(\mathbf x)=\delta_k(r),k=1,\ldots,m,$ the following asymptotic equality holds 
\begin{multline}\label{eq:4}
\delta_k(r)=\delta_k(r)\Psi_k^{n-2}(r)(1+\Psi_k^{n-2}(r)R_k^{2-n})(|\mathbf x-\mathbf x_k|^{2-n}-R_k^{2-n}+o(1))\\
-\delta_k(r)\sum \limits_{{l=1,} \atop{l\not=k}}^m\Psi_k^{n-2}(r)\Psi_l^{n-2}(r)g(\mathbf x_k,\mathbf x_l)(g(\mathbf x_k,\mathbf x_l)+o(1))\\
-\sum \limits_{{j=1,} \atop {j\not=k}}^m\delta_j(r)\Psi_k^{n-2}(r)\Psi_j^{n-2}(r)g(\mathbf x_k,\mathbf x_j)(|\mathbf x-\mathbf x_k|^{2-n}-R_k^{2-n}+o(1))\\
+\sum \limits_{{j=1,} \atop {j\not=k}}^m\delta_j(r)\Psi_j^{n-2}(r)(1+\Psi_j^{n-2}(r)R_j^{2-n})(g(\mathbf x_k,\mathbf x_j)+o(1))\\
-\sum \limits_{{j=1,} \atop {j\not=k}}^m\delta_j(r)\sum \limits_{{l=1,l\not=j,} \atop {l\not=k}}^m\Psi_l^{n-2}(r)\Psi_j^{n-2}(r)g(\mathbf x_l,\mathbf x_j)(g(\mathbf x_k,\mathbf x_l)+o(1)),\mathbf x\to \mathbf x_k.
\end{multline}
Since the product $\Psi_k^{n-2}(r)|\mathbf x-\mathbf x_k|^{2-n}$  is bounded, we conclude that 
\begin{equation}\label{eq:add}
\delta_k(r)=\delta_k(r)\Psi_k^{n-2}(r)(1+\Psi_k^{n-2}(r)R_k^{2-n})(|\mathbf x-\mathbf x_k|^{2-n}-R_k^{2-n})+o(1).
\end{equation}
It follows that 
\begin{multline}
|\mathbf x-\mathbf x_k|^{2-n}\Psi_k^{n-2}(r)-R_k^{2-n}\Psi_k^{n-2}(r)+o(1)=(1+\Psi_k^{n-2}(r)R_k^{2-n})^{-1}(1+o(1))\\=1-R_k^{2-n}\Psi_k^{n-2}(r)+o(1).
\end{multline}
 Substituting the third sum term \eqref{eq:4} for the equality $|\mathbf x-\mathbf x_k|^{2-n}\Psi_k^{n-2}(r)=1+o(1)$, we obtain 
 \begin{multline}\label{eq:6}
\delta_k(r)=\delta_k(r)\Psi_k^{n-2}(r)(1+\Psi_k^{n-2}(r)R_k^{2-n})(|\mathbf x-\mathbf x_k|^{2-n}-R_k^{2-n})\\-\sum \limits_{{j=1,} \atop {j\not=k}}^n\delta_j(r)\Psi_j^{n-2}(r)g(\mathbf x_k,\mathbf x_j)(1-\Psi_k^{n-2}(r)R_k^{2-n}+o(1)) \\ 
+\sum \limits_{{j=1,} \atop {j\not=k}}^n\delta_j(r)\Psi_j^{n-2}(r) g(\mathbf x_k,\mathbf x_j)+o(r^{n-2}).
\end{multline}
Rearranging \eqref{eq:6}, we obtain 
\begin{multline}
|\mathbf x-\mathbf x_k|^{2-n}\Psi_k^{n-2}(r)-R_k^{2-n}\Psi_k^{n-2}(r)+o(r^{n-2})=(1+\Psi_k^{n-2}(r)R_k^{2-n})^{-1}\\=1-R_k^{2-n}\Psi_k^{n-2}(r)+o(r^{n-2}).
\end{multline}
Therefore, in the points $\mathbf x\in \partial E_k,$ where $g(\mathbf x)=\delta_k(r),$ the asymptotic equality holds  
$$ \frac{|\mathbf x-\mathbf x_k|}{\mu_k r}=\frac{|\mathbf x-\mathbf x_k|}{\Psi_k(r)}=(1+o(r^{n-2}))^{1/(2-n)}=1+O(r^{n-1}),\ r \to 0.$$ Taking into account the harmonicity of $g(\mathbf x)$, we conclude that  at the boundary of $\tilde E_k(r)=\{\mathbf x\in B: g(\mathbf x)/\delta_k(r)\geq 1\}$ we have $g(\mathbf x)=\delta_k(r),$ $ |\mathbf x-\mathbf x_k|=\mu_k r(1+O(r^{n-1})),$
($r$ is sufficiently small). 

Let $B_r=B\textbackslash\bigcup\limits_{k=1}^n \tilde E_k(r)$, $\tilde C (r)=(\{\Gamma,\tilde E _1(r),\ldots,\tilde E_m(r)\}, \{0, \delta_1(r),\ldots, \delta_m(r)\}, B\}).$ Using the Green's formula and the Gauss' theorem, we have 
\begin{multline}
cap\,\tilde C(r) = \int\limits_B |\nabla g|^2d\mathbf x =-\int\limits_{\partial B_r} g\frac{\partial g}{\partial n}d\sigma = -\sum\limits_{k=1}^m\delta_k(r)\int\limits_{\partial B_r} \frac{\partial g}{\partial n}d\sigma\\=-\sum\limits_{k=1}^m\delta_k(r)\sum\limits_{l=1}^m\delta_l(r)\beta_{kl}\int\limits_{|\mathbf x-\mathbf x_k|=\rho}\frac{\partial |\mathbf x-\mathbf x_k|^{2-n}}{\partial n}d\sigma = -\sum\limits_{k=1}^m\delta_k(r)\sum\limits_{l=1}^m\delta_l(r)\beta_{kl}(2-n)\omega_{n-1}\\=
(n-2)\omega_{n-1}\sum\limits_{k=1}^m\delta_k(r)\sum\limits_{l=1}^m\delta_l(r)\beta_{kl},
\end{multline}
where $n$ is the inward normal.
Substituting with $\beta_{kl}$ values and considering the fact that condenser capacity does not depend on $\rho$, we obtain
\begin{multline}\label{eq:10}
\frac{1}{(n-2)\omega_{n-1}}cap\,\tilde C(r) = 
\sum\limits_{k=1}^m\delta_k^2(r)\Psi_k^{n-2}(r)(1+\Psi_k^{n-2}(r)R_k^{2-n})\\ - \sum\limits_{k=1}^m\sum \limits_{{l=1,} \atop{l\not=k}}^m\delta_k(r)\delta_l(r)\Psi_k^{n-2}(r)\Psi_l^{n-2}(r)g(\mathbf x_k,\mathbf x_l).
\end{multline}
Considering $\delta_k(r)=\delta_k+\alpha_kr^{n-2}+o(r^{n-2})$, we obtain 
\begin{multline}\label{eq:14}
\frac{1}{(n-2)\omega_{n-1}}cap\,\tilde C(r) = 
\sum\limits_{k=1}^m(\delta_k^2+2\alpha_k\delta_kr^{n-2}+o(r^{n-2}))(1+(\mu_kr)^{n-2}R_k^{2-n})\mu_k^{n-2}r^{n-2}\\-\sum\limits_{k=1}^m\sum \limits_{{l=1,} \atop{l\not=k}}^m[\delta_k\delta_l+(\alpha_k\delta_l+\alpha_l\delta_k)r^{n-2}+o(r^{n-2})](\mu_kr)^{n-2}(\mu_lr)^{n-2}g(\mathbf x_k,\mathbf x_l)=
\left(\sum\limits_{k=1}^m\delta_k^2\mu_k^{n-2}\right)r^{n-2}\\+\left(\sum\limits_{k=1}^m(\delta_k^2\mu_k^{2n-4}R_k^{2-n}+2\alpha_k\delta_k\mu_k^{n-2})-\sum\limits_{k=1}^m\sum\limits_{{l=1,} \atop{l\not=k}}^m \delta_k\delta_l\mu_k^{n-2}\mu_l^{n-2}g(\mathbf x_k,\mathbf x_l)\right)r^{2n-4}+o(r^{2n-4}).
\end{multline}

Let $r^\pm=r(1\pm\sqrt{r}r^{n-2})=r(1+o(r^{n-2})),$  $\delta_k^\pm(r)=\delta_k(r)(1\pm\sqrt{r}r^{n-2})=\delta_k(1+o(r^{n-2})).$ Let us construct $\tilde E_k^\pm(r)$ sets as in the case with $\tilde E_k(r)$ repeating our previous reasoning by substituting $r$ and $\delta_k(r)$ with $r^+$ and $\delta^+_k(r).$  At the $\tilde E_k^\pm(r)$ boundary it holds that 
$$ \frac{|\mathbf x-\mathbf x_k|}{\mu_k r}=(1\pm\sqrt{r}r^{n-2})(1+O(r^{n-1})).$$ 
That is why the following relations hold assuming $r$ is sufficiently small:
$$\tilde E_k^-(r)\subset F_k\subset \tilde E_k^+(r),\ 
\delta^+_k(r)/\delta_k(r)\geq 1,\ \delta_k^-(r)/\delta_k(r)\leq 1.$$ Therefore, for the condenser capacities
$$C^\pm (r)=(\{\Gamma,\tilde E^\pm_1(r),\ldots,\tilde E_m^\pm(r)\}, \{0, \delta_1^\pm(r),\ldots, \delta_m^\pm(r)\}, B\})$$ and 
$$C(r)=(\{\Gamma, F_1,\ldots,F_m\}, \{0, \delta_1 (r),\ldots, \delta_m(r)\}, B\})$$ the following inequalities hold 
\begin{equation}\label{eq:13}
cap\, C^-(r)\leq cap\, C(r)\leq cap\, C^+(r).
\end{equation}
From the formula \eqref{eq:10}, it follows that
\begin{multline}\label{eq:12}
\frac{1}{(n-2)\omega_{n-1}}cap\,C^\pm(r) = \\
\sum\limits_{k=1}^m\delta_k^\pm(r)^2\Psi_k^{n-2}(r^\pm)(1+\Psi_k^{n-2}(r^\pm)R_k^{2-n})- \sum\limits_{k=1}^m\sum \limits_{{l=1,} \atop{l\not=k}}^m\delta_k^\pm(r)\delta_l^\pm(r)\Psi_k^{n-2}(r^\pm)\Psi_l^{n-2}(r^\pm)g(\mathbf x_k,\mathbf x_l)=\\=
\sum\limits_{k=1}^m\delta_k(r)^2\Psi_k^{n-2}(r)(1+\Psi_k^{n-2}(r)R_k^{2-n}+o(r^{n-2})) -\\ \sum\limits_{k=1}^m\sum \limits_{{l=1,} \atop{l\not=k}}^m\delta_k(r)\delta_l(r)\Psi_k^{n-2}(r)\Psi_l^{n-2}(r)g(\mathbf x_k,\mathbf x_l)(1+o(r^{n-2}))+o(r^{2n-4})=\\= 
\sum\limits_{k=1}^m\delta_k(r)^2\Psi_k^{n-2}(r)(1+\Psi_k^{n-2}(r)R_k^{2-n})- \sum\limits_{k=1}^m\sum \limits_{{l=1,} \atop{l\not=k}}^m\delta_k(r)\delta_l(r)\Psi_k^{n-2}(r)\Psi_l^{n-2}(r)g(\mathbf x_k,\mathbf x_l)+\\o(r^{2n-4})=\frac{1}{(n-2)\omega_{n-1}}cap\,\tilde C(r)+o(r^{2n-4}).
\end{multline} Considering \eqref{eq:13}, $cap\, C(r)=cap\,\tilde C(r)+o(r^{2n-4}).$ Now the assertion of the theorem follows from \eqref{eq:14}. QED.

Considering Theorem 1 in the case $B=B(0,R)=\{|\mathbf x|<R\}$, $\Gamma=\partial B(0,R)$, $X=\{0\},$ $\Delta=\{\delta\}$, $\Omega=\{\alpha\}$, $\Psi=\{1\}$, the reduced modulus equals $$M(B(0,R),\partial B(0,R),\{0\},\{\delta\},\{\alpha\},\{1\})=-\frac{\lambda_n}{\delta^2}R^{2-n}-\frac{2\alpha\lambda_n}{\delta^3}.$$ In this case it is easy to check the obtained result by calculating the limit defined in \eqref{eq:mod2}. In effect, according to the well-known formula for the calculation of condenser capacity with sphere plates \cite[p.97]{Maz} 
$$|C(r;B(0,R),\partial B(0,R),\{0\},\{\delta\},\{\alpha\},\{1\})|=\frac{\lambda_n(r^{2-n}-R^{2-n})}{(\delta+\alpha r^{n-2}+o(r^{n-2}))^2}.$$ Therefore, $$\lim\limits_{r \to 0}(|C(r;B(0,R),\partial B(0,R),\{0\},\{\delta\},\{\alpha\},\{1\})|-\delta^{-2}\lambda_n r^{2-n})=-\frac{\lambda_n}{\delta^2}R^{2-n}-\frac{2\alpha\lambda_n}{\delta^3}.$$

In a series of questions it turns out to be useful to represent the reduced modulus directly in the terms of the Dirichlet's integral of the so-called potential function \cite{DubPril5}. As the following theorem demonstrates, now \eqref{eq:potential} plays the role of the potential function for the modulus $M(B,\Gamma,X,\Delta,\Omega,\Psi)$.

\begin{theorem}\label{Th2}
Considering Theorem 1, 
\begin{equation}
M(B,\Gamma,X,\Delta,\Omega,\Psi)=\lim_{r \to 0}\left( \nu^2 I(u,B_r) -\nu\lambda_n r^{2-n}\right),
 \end{equation}
where $B_r=B\setminus\bigcup_{l=1}^m E(x_l,\mu_l r),$  \begin{equation}\label{eq:potential}
u(\mathbf x) = \sum\limits_{l=1}^m \mu^{n-2}_l (\delta_l-\alpha_l r^{n-2})g_B(\mathbf x,\mathbf x_l,\Gamma).
\end{equation}
\end{theorem} 
\textbf{Proof.} From the Green's identity
\begin{equation}
\int_V|\nabla u|^2d\mathbf x = -\int_{\partial V}  u\frac{\partial u}{\partial n} d\sigma,
\end{equation}
it follows that
\begin{equation}\label{eq:IntU}
I(u,B_r) =-\int_{\partial B_r}  u\frac{\partial u}{\partial n} d\sigma=-\sum\limits_{k=1}^m \int_{|\mathbf x-\mathbf x_k|=\rho_k}  u\frac{\partial u}{\partial n} d\sigma,
\end{equation}
where $\rho_k=\mu_k r.$ The second equation in \eqref{eq:IntU} holds, because $u\frac{\partial u}{\partial n} =0$ at ${\partial B}$. Let $\tilde \delta_k(r)=\delta_k-\alpha_k r^{n-2}.$ 
Considering $|\mathbf x-\mathbf x_k| = \rho_k$,  it holds that 
\begin{equation}
u(\mathbf x) = \mu^{n-2}_k \tilde\delta_k(r)\lambda_n\rho_k^{2-n}+h_k(\mathbf x),
\end{equation}
where 
$h_k(\mathbf x)$ are functions harmonic in the neighborhoods of $\mathbf x_k$.
\begin{multline}\label{eq:int1}
I(u,B_r) =-\sum\limits_{k=1}^m \int_{S(\mathbf x_k,\rho_k)}  u\frac{\partial u}{\partial n} ds=\\-\sum\limits_{k=1}^m \int_{S(\mathbf x_k,\rho_k)}
\left(\mu^{n-2}_k \tilde\delta_k(r)\lambda_n\rho_k^{2-n}+h_k(\mathbf x)\right)\left(-\frac{\mu_k^{n-2}\tilde\delta_k(r)\rho_k^{1-n}}{w_{n-1}}+\frac{\partial h_k(\mathbf x)}{\partial n} \right)d\sigma=\\\sum\limits_{k=1}^m \int_{S(\mathbf 0,1)}
\left(\mu^{n-2}_k \tilde\delta_k(r)\lambda_n\rho_k^{2-n}+h_k(\mathbf x)\right)\left(\frac{\mu_k^{n-2}\tilde\delta_k(r)}{w_{n-1}}-\rho_k^{n-1}\frac{\partial h_k(\mathbf x)}{\partial n} \right)d\sigma.
\end{multline}
 Assuming in \eqref{eq:int1} $\rho_k=\mu_k r$ we obtain 
\begin{equation}\label{IntDir} 
I(u,B_r) = r^{2-n}\lambda_n\sum\limits_{k=1}^m \mu^{n-2}_k \tilde \delta^2_k(r)+\sum\limits_{k=1}^m \mu^{n-2}_k \tilde \delta_k(r)a_k(r)+o(1), r\to 0, 
\end{equation}
where 
\begin{equation}\label{ak}
a_k(r)=-\mu^{n-2}_k \tilde\delta_k(r)\lambda_n r(B,\Gamma,\mathbf x_k)^{2-n}+\sum \limits_{{l=1,} \atop{l\not=k}}^m\mu^{n-2}_l \tilde\delta_l(r)g_B(\mathbf x_k,\mathbf x_l,\Gamma),
\end{equation}
From the asymptotics of $\tilde\delta_k(r)$ and the expansion \eqref{ak}, it follows that
$$  r^{2-n}\lambda_n\sum\limits_{k=1}^m \mu^{n-2}_k \tilde \delta^2_k(r)=\frac{r^{2-n}\lambda_n}{\nu}-2\sum_{k=1}^m\alpha_k\delta_k\mu_k^{n-2}\lambda_n+o(1),$$ $$
\sum\limits_{k=1}^m \mu^{n-2}_k \tilde \delta_k(r)a_k(r)=\frac{M(B,\Gamma,X,\Delta,\Omega,\Psi)}{\nu^2}+2\sum_{k=1}^m\alpha_k\delta_k\mu_k^{n-2}\lambda_n+o(1), r\to 0.$$
Now the assertion of the theorem follows from \eqref{IntDir}. QED.

\section{The relation between the reduced modulus and the discrete energy. The behavior during some geometric transformations}

Let $\Delta=\{\delta_k\}_{k=1}^m$ be a discrete charge
(the set of real numbers other than zero), which has the value of $\delta_k$ in the point $\mathbf x_k\in X$, where $k=1,\ldots,m$. We call the \emph{Robin energy} of this charge with respect to the domain $B$ and the set $\Gamma$ the following quantity 

$$
E(X,\Delta,B,\Gamma)=\sum\limits_{k=1}^m\sum \limits_{{l=1} \atop {l\not=
k}}^m \delta_k\delta_l g_{\small B}(\mathbf x_k,\mathbf x_l,\Gamma).
$$ 
If  $\Gamma=\partial B$ the Robin energy coincides with the Green energy \cite{Dubinin1}. 
 \begin{theorem} 
Considering Theorem 1, 
\begin{equation}
E(X,\Delta,B,\Gamma)=(\sum_{k=1}^m\delta _k^2)^2M(B,\Gamma,X, \Delta,\{-\delta_kr(B,\Gamma,\mathbf x_k)^{2-n}/2\}_{k=1}^m,\{1\}_{k=1}^m),
\end{equation}
\begin{equation}\label{Energy}
E(X,\Delta,B,\Gamma)=(\sum_{k=1}^m\delta _k^2)^2M(B,\Gamma,X, \Delta,\{0\}_{k=1}^m,\{1\}_{k=1}^m)+\sum_{k=1}^m\delta_k^2 \lambda_n r(B,\Gamma,\mathbf x_k)^{2-n}.
\end{equation}
\end{theorem}
\textbf{Proof} follows directly from the formula \eqref{basic}.

Let  $D$ and $B$ be domains, such that $B\subset D$ and $(\partial B)\cap D \subset\Gamma.$ In this case we say that the domain $D$ is obtained from the domain $B$ expanding the boundary part $\Gamma\subset \partial B.$
\begin{theorem}  Let $B_1\subset B_2$, and the domain $B_2$ is obtained from the    domain $B_1$ expanding the boundary part $\Gamma_1\subset \partial B_1,$
and the set $\Gamma_2\subset \partial B_2$ satisfies the condition $\Gamma_2\subset( \Gamma_1\cup (\mathbb{R}^n\setminus \overline{B_1})),$ $X\subset B_1.$
Then 
\begin{equation}\label{eq:innn}M(B_1,\Gamma_1,X,\Delta,\Omega,\Psi)\leq M(B_2,\Gamma_2,X,\Delta,\Omega,\Psi).\end{equation}
If the domain $B_2$ is obtained from the domain $B_1$ expanding the boundary part $\partial B_1\setminus \Gamma_1$ (complement $\Gamma_1$ to the boundary $\partial B_1$), and $\Gamma_1\subset \Gamma_2,$ then \begin{equation}\label{eq:in2}M(B_1,\Gamma_1,X,\Delta,\Omega,\Psi)\geq M(B_2,\Gamma_2,X,\Delta,\Omega,\Psi).\end{equation}
\end{theorem}

\textbf{Proof}.  Let $v_1$ be  admissible  for the condenser capacity $C(r;B_1,\Gamma_1,$ $ X,\Delta,\Omega,\Psi)$. We  define $v_1$ as zero-equal everywhere else except for $B_1.$ So, $v_1$ is  admissible for  $C(r;B_2,\Gamma_2, X,\Delta,$ $ \Omega,\Psi)$ in the case of expanding $B_1$ by the boundary part $\Gamma_1.$ In the case of expanding complement $\Gamma_1$ to the boundary $\partial B_1$, we are considering the admissible function $v_2$ for  $C(r;B_2,\Gamma_2, X,\Delta,$ $ \Omega,\Psi).$ We conclude that $v_2$ is  admissible  for the $C(r;B_1,\Gamma_1, X,\Delta,\Omega,\Psi)$. Henceforth, inequalities \eqref{eq:innn}, \eqref{eq:in2} follow from the definitions of the capacity and the  reduced modulus. QED.                                                      

For "ordinary"\ condensers  with two plates a series of symmetrization transformations that do not increase capacity is known. The symmetrization transformations of the generalized condensers with several plates and various potential levels are barely studied. For some potential choices and plate geometries the proof of the symmetrization principle is extended to the case of several plates without many adjustments. An example we can give here is the dissymmetrization that we have used in the proof of Theorem 1 in \cite{DubPrilGreen}  when estimating  the Green energy of the rotation body. Let us give a similar result related to the Robin energy. We will need cylindrical coordinates $(r,\theta, \mathbf{x}')$ of a point
$\mathbf{x}=(x_1,\dots,\,x_n)$ in $\mathbb{R}^n$, which are related to the Cartesian coordinates by $x_1=r \cos \theta,$ $x_2=r \sin \theta,$
$\mathbf{x}'\in J$,    $J$ -- $(n-2)$--dimensional plane $\{\mathbf
x\in\mathbb{R}^n:\mathbf x=(0,0,x_3,\ldots, x_n)\}$. The designation $\{\theta=\phi\}$ denotes the set of points in $\mathbb{R}^n$ with coordinates $(r,\phi, \mathbf{x}'),r\geq 0,\mathbf{x}'\in J$, where $\phi$ is fixed. We will call set $Q$ the rotation body relative to $J$ if for all $\varphi$ not only  $(r,\theta, \mathbf{x}')\in B$, but also  $(r,\theta+\varphi, \mathbf{x}')\in B.$ For instance, rotation bodies  are spheres $\{|\mathbf{x}|<\tau\}$, rings $\{|\tau_1<|\mathbf{x}|<\tau_2\}$, hyperspheres $\{|\mathbf{x}|=\tau\}$, etc. 

\begin{theorem}
Let $B$  and  $\Gamma$ be  rotation bodies, $\Omega=\{S\}$ be the set comprised of a finite number of different circles   $S=\{(r_0,\theta,\mathbf
x'_0):0\leq\theta\leq 2\pi\}$, which lie in $B$ (here
$r_0>0$ and $\mathbf \mathbf{x}'_0\in J$ are fixed). For arbitrary real numbers $\theta_j,$ $j=0,\ldots,$ \makebox{$m-1,$}
\begin{equation*}
0\leq \theta_0<\theta_1<\ldots<\theta_{m-1}<2\pi,
\end{equation*}
let us denote by $X=\{\mathbf x_k\}_{k=1}^p$ -- the set of the intersection points of the circles  from $\Omega$ with the half-planes
\begin{equation*}
L_j=\{(r,\theta,\mathbf x'):\theta=\theta_j\}, \ j=0,\ldots,m-1,
\end{equation*}
  $X^*=\{\mathbf x^*_k\}_{k=1}^p$ -- the set of the intersection points of the circles  from  $\Omega$ with the symmetrical  half-planes
\begin{equation*}
\ L_j^*=\{(r,\theta,\mathbf x'): \theta=2\pi j/m\}, \
j=0,\ldots,m-1.
\end{equation*}
 Let the charge $\Delta=\{\delta_k\}_{k=1}^p$ have the same values $\delta_k=\delta_l$ in the points $\mathbf x_k\in X$ and
$\mathbf x_l\in X$, which are located on the same circle from $\Omega$ and, additionally, let the points $\mathbf x_k\in X$
and $\mathbf x^*_k\in X^*$ be located on the same circle in $\Omega$, where
$k=1,\ldots,p$. Then
$$
E(X,\Delta,B,\Gamma)\geq E(X^*,\Delta,B, \Gamma).
$$
\end{theorem}
\textbf{Proof.} In the proof of Theorem 1 from \cite{DubPrilGreen} let us substitute the condenser $C(t,B,$ $ X^*,\Lambda, \Psi)$ with the condenser $C(r;B,\Gamma, X^*,\Delta,\{0\}_{k=1}^m,\{1\}_{k=1}^m)$, and the asymptotic formula (5) from  \cite{DubPrilGreen} -- with the formula \eqref{ascon}. QED. 

Speaking about the proof of the monotonicity of the reduced modulus during geometric transformations, there exists an additional opportunity to use the relation between the reduced modulus and the potential function. Then we can obtain the results for the reduced moduli in the cases where it is impossible to prove the symmetrization principle for the capacities of condensers of three or more plates. The idea of this approach traces back to  \cite{DubPril5}. In this paper we will construct a polarization. Let us start with its definition.

Let the hyperplane $L$ break the space $\mathbb R^n$ into two open non-intersecting half-spaces $\mathbf R^+$, $\mathbf R^-$, such that $\mathbb R^n=\mathbf R^+\cup \mathbf R^-\cup L.$ Let $\mathbf x^*$ denote the point symmetrical to $\mathbf x$ relative to the hyperplane $L.$ 
The polarization of the set $Q\subset \mathbb R^n$ is defined as the transition to the set  $PQ=(Q\cup Q^*)^+\cup(Q\cap Q^*)^-\cup (L\cap Q).$ Here $Q^*$ denotes the set symmetrical to $Q$ relative to $L,$  $Q^{\pm}=Q\cap\mathbf R^\pm.$

 For the collections $X=\{\mathbf x_k\}_{k=1}^m$, $\Delta=\{\delta_k\}_{k=1}^m$,  $\Omega=\{\alpha_k\}_{k=1}^m$, $ \Psi=\{\mu_k\}_{k=1}^m$ let us define the polarization relative to $L$ as a transition to the sets  $PX=\{P\mathbf x_k\}_{k=1}^m$, $P\Delta=\{P\delta_k\}_{k=1}^m$,  $P\Omega=\{P\alpha_k\}_{k=1}^m$, $ P\Psi=\{P\mu_k\}_{k=1}^m$ according to the following rules. We will allow to construct the polarization only under the additional condition:  if there are symmetric points $x_k$, $x_l$ in the set $X$, $x_k\neq x_l,$ then $\mu_k^{n-2}\delta_k\neq\mu_l^{n-2}\delta_l$. Let $B=B^+\cup B^-\cup (L\cap B).$ In the case $\mathbf x_k\in L$ four-tuple $(\mathbf x_k,\delta_k,\alpha_k,\mu_k)$ keeps unchanged, i.e.  $(P\mathbf x_k,P\delta_k,P\alpha_k,P\mu_k)=(\mathbf x_k,\delta_k,\alpha_k,\mu_k).$ If $\mathbf x_k\in B^+,$ $\mathbf x_k^*\notin X$, then the four-tuple $(\mathbf x_k,\delta_k,\alpha_k,\mu_k)$ keeps unchanged in the case of a positive charge $\delta_k>0$, and in the case of a negative charge $\delta_k<0$ the transition to the symmetrical point takes place: $(P\mathbf x_k,P\delta_k,P\alpha_k,P\mu_k)=(\mathbf x_k^*,\delta_k,$ $\alpha_k,\mu_k).$ If $\mathbf x_k\in B^-,$ $\mathbf x_k^*\notin X$, then the four-tuple $(\mathbf x_k,\delta_k,\alpha_k,\mu_k)$ keeps unchanged in the case of a negative charge $\delta_k<0$, and in the case of a positive charge the transition to the symmetrical point takes place $\delta_k>0$. If $\mathbf x_k^*=\mathbf x_l\in X,$ $\mathbf x_k\notin L,$ then the polarization of a pair of four-tuples  $(\mathbf x_k,\delta_k,\alpha_k,\mu_k),$ $(\mathbf x_l,\delta_l,\alpha_l,\mu_l)$ will be jointly defined. Let $P\mathbf x_k$ denote the point of the pair $\mathbf x_k,$ $\mathbf x_l$ from $B^+,$ and let $P\mathbf x_l$ denote the point from $B^-.$
Then, let $(P\delta_k,P\alpha_k,P\mu_k)$ denote one of the three-tuples $(\delta_k,\alpha_k,\mu_k)$ and $(\delta_l,\alpha_l,\mu_l)$, for which  $\mu_k^{n-2}\delta_k> \mu_l^{n-2}\delta_l,$   and let us denote $(P\delta_l,P\alpha_l,P\mu_l)$ the remaining three-tuple.  
 \begin{theorem} 
Let $B$ and $\Gamma$ be symmetrical relative to  $L.$ Then \begin{equation}\label{eq:inM}M(B,\Gamma,X,\Delta,\Omega,\Psi)\leq M(B,\Gamma,PX,P\Delta,P\Omega,P\Psi),\end{equation}
\begin{equation}\label{eq:in1}E(X,\Delta,B,\Gamma)\leq E(PX,P\Delta ,B,\Gamma),\end{equation}
\end{theorem}

 {\bf{Proof.}} Let us denote $ \tilde{\delta}_l(r)=\delta_l- \alpha_l r^{n-2},$  $ {\widetilde{P\delta}}_l(r)=P\delta_l-(P\alpha_l) r^{n-2}$, where $l=1,\ldots m.$ Let us examine the following potential functions $$
u(\mathbf x) = \sum\limits_{l=1}^m \mu^{n-2}_l \tilde\delta_l(r)g_B(\mathbf x,\mathbf x_l,\Gamma),
\ \ \ u_1(\mathbf x) = \sum\limits_{l=1}^m (P\mu_l)^{n-2} {\widetilde{ P\delta}}_l(r) g_B(P\mathbf x,P\mathbf x_l,\Gamma),
$$
of the reduced moduli $M(B,\Gamma,X,\Delta,\Omega,\Psi),$ $M(B,\Gamma,PX,P\Delta,P\Omega,P\Psi)$, respectively. As it is shown in the proof of the theorem \ref{Th2}, the following asymptotic expansions at $r\to 0$  hold 
 \begin{equation}
u(\mathbf x) = \tilde\delta_k(r)\lambda_n r^{2-n}+a_k(r)+o(1),\  |\mathbf x-\mathbf x_k|=\mu_k r,
\end{equation}
\begin{equation}
u_1(\mathbf x) = \widetilde{P\delta}_k(r)\lambda_n r^{2-n}+b_k(r)+o(1),  \ |\mathbf x-P\mathbf x_k|=(P\mu_k) r,
\end{equation}
where  
\begin{equation}\label{sumcon1}
\sum\limits_{k=1}^m \mu^{n-2}_k \tilde \delta_k(r)a_k(r)=\frac{M(B,\Gamma,X,\Delta,\Omega,\Psi)}{\nu^2}+2\sum_{k=1}^m\alpha_k\delta_k\mu_k^{n-2}\lambda_n+o(1), r\to 0.
\end{equation}
\begin{equation}\label{sumcon2}
\sum\limits_{k=1}^m (P\mu_k)^{n-2} \widetilde{P\delta_k}(r)b_k(r)=\frac{M(B,\Gamma,PX,P\Delta,P\Omega,P\Psi)}{\nu^2}+2\sum_{k=1}^m (P\alpha_k) (P\delta_k) (P\mu_k)^{n-2}\lambda_n+o(1).\end{equation}

Let us examine the function  
\begin{equation}
v(\mathbf x) = \left\{\begin{array}{l}
\max (u(x),u(x^*)),\ x\in B^+\cup(B\cap L)\\
\min{(u(x),u(x^*)}),\ x\in B^-,
\end{array}\right.
\end{equation}
 As $g_B(\mathbf x,\mathbf x_k, \Gamma)\to +\infty$ assuming $\mathbf x\to\mathbf x_k,$ the definition of the polarization of the collections $X$,  $\Delta$,  $\Omega$, $ \Psi$ leads to the following asymptotics  $r\to 0$ \begin{equation}
v(\mathbf x) ={\widetilde {P\delta}}_k(r)\lambda_n r^{2-n}+c_k(r)+o(1),  \ |\mathbf x-P\mathbf x_k|=(P\mu_k)r,
\end{equation} 
where $$
\sum\limits_{k=1}^m (P\mu_k)^{n-2} (\widetilde {P\delta_k}(r))c_k(r)=\sum\limits_{k=1}^m \mu^{n-2}_k \tilde \delta_k(r)a_k(r).$$

Let us denote $B_r=B\setminus\bigcup_{l=1}^m E(x_l,\mu_l r),$ $PB_r=PB\setminus\bigcup_{l=1}^m E(Px_l,(P\mu_l)r).$ From the definition of the function $v$, it follows that 
$$I(v,PB_r)=I(u, B_r).$$ Using the Green's formula and the fact that $ \int_{\partial B}(v-u_1)\frac{\partial u_1}{\partial n} d\sigma=0,$ we obtain  

\begin{equation}\label{eq:Posit}
0\leq I(v-u_1,PB_r)=I(v,PB_r)-I(u_1,PB_r)+2\sum_{k=1}^m\int_{S(P\mathbf x_k,(P\mu_k)r)}(v-u_1)\frac{\partial u_1}{\partial n}d\sigma.
\end{equation}
Repeating the proofs of the inequalities   \eqref{IntDir}, we obtain  
\begin{multline}
I(v,PB_r)-I(u_1,PB_r)=I(u,B_r)-I(u_1,PB_r) =\\
\sum\limits_{k=1}^m \mu^{n-2}_k \widetilde{\delta_k}(r)a_k(r)-\sum\limits_{k=1}^m P\mu^{n-2}_k \widetilde{P\delta_k}(r)b_k(r)+o(1), \ r\to 0,
\end{multline}
\begin{multline}
\sum_{k=1}^m\int_{S(P\mathbf x_k,(P\mu_k)r)}(v-u_1)\frac{\partial u_1}{\partial n}d\sigma=\sum\limits_{k=1}^m (P\mu_k)^{n-2} \widetilde{P\delta_k}(r)(b_k(r)-c_k(r))+o(1)=\\\sum\limits_{k=1}^m (P\mu_k)^{n-2}\widetilde{ P\delta_k}(r)b_k(r) -\sum\limits_{k=1}^m \mu^{n-2}_k \tilde \delta_k(r)a_k(r)+o(1), \ r\to 0.
\end{multline}
From  \eqref{eq:Posit} we obtain 
$$0\leq \sum\limits_{k=1}^m (P\mu_k)^{n-2}\widetilde{ P\delta_k}(r)b_k(r) -\sum\limits_{k=1}^m \mu^{n-2}_k \tilde \delta_k(r)a_k(r)+o(1), \ r\to 0.$$ Taking into account \eqref{sumcon1}, \eqref{sumcon2} we obtain \eqref{eq:inM}.  The inequality \eqref{eq:in1} follows from the formula \eqref{Energy} and the fact that due to $B$ and $\Gamma$ symmetry 
$$\sum_{k=1}^m\delta_k^2 r(B,\Gamma,\mathbf x_k)^{2-n}=\sum_{k=1}^m (P\delta_k)^2 r(B,\Gamma,P\mathbf x_k)^{2-n}.$$
QED. 

 \begin{theorem} 
If $\delta_k>0,$ $k=1,\ldots,m,$ and $PB$ is a domain, then \begin{equation}M(B,\partial B,X,\Delta,\Omega,\Psi)\leq M(PB,\partial (PB),PX,P\Delta,P\Omega,P\Psi).\end{equation}
\end{theorem}

\textbf{Proof.} The proof is almost identical to the one of the previous theorem with a single remark that assuming the charges are positive and $\Gamma=\partial B$, the function $u$ is non-negative in $B$, and it can be extended as null-valued outside $B$ due to its continuity.  

{\bf{Acknowledgments.}}
The research is supported by the RSF (project No. 23-21-00056, https://rscf.ru/project/23-21-00056/)

\section*{Statements and Declarations}

{\bf Data Availability Statement:} Not applicable.\\
{\bf Conflicts of Interest:} The authors declare no conflicts of interest.

%
%

\end{document}